\documentclass[11pt,english,oneside]{amsart}
\usepackage[T1]{fontenc}
\usepackage[latin9]{inputenc}
\usepackage{geometry}
\usepackage{amssymb}
\usepackage{color}

\usepackage{graphicx,color}
\usepackage{amsmath, amssymb, graphics}

\usepackage{graphicx}

\makeatletter
\numberwithin{equation}{section} 
\numberwithin{figure}{section} 
  \@ifundefined{theoremstyle}{\usepackage{amsthm}}{}
  \theoremstyle{plain}
  \newtheorem{thm}{Theorem}[section]
  \theoremstyle{plain}
  \newtheorem{cor}[thm]{Corollary}
  \theoremstyle{plain}
  \newtheorem{prop}[thm]{Proposition}
  \theoremstyle{remark}
  \newtheorem{rem}[thm]{Remark}
  \theoremstyle{remark}
  
  \theoremstyle{plain}
  \newtheorem{lem}[thm]{Lemma}

\usepackage{geometry}

\def\bfR#1{{\bf R}^#1}

\def\com#1{ \hbox{#1}}

\def\eop{{\vrule height 6pt width 5pt depth 0pt}}\smallskip
\def\<{{\langle }}
\def\>{{\rangle }}

\makeatother

\usepackage{babel}

\def\bfR#1{{\bf R}^#1}

\def\com#1{ \quad\hbox{#1}\quad}

\def\eop{{\vrule height 6pt width 5pt depth 0pt}}\smallskip
\def\<{{\langle }}
\def\>{{\rangle }}

\makeatother

\begin{document}

\title{A dynamical interpretation of the profile curve of cmc Twizzlers surfaces}

\author{ Oscar M. Perdomo }

\date{\today}

\curraddr{Department of Mathematics\\
Central Connecticut State University\\
New Britain, CT 06050\\
}

\email{ perdomoosm@ccsu.edu}

\begin{abstract}

It is known that for any non-zero $M\in (-\frac{1}{4},\infty)$, if we roll the conic $\{(x,y): 4 x^2-\frac{y^2}{M}=1\}$ on a line in a plane, and then we rotate about this line the trace of a focus, then we obtain a surface of revolution $\mathbb{D}(M)$ with mean curvature 1. If $M<0$, $\mathbb{D}(M)$ is embedded and it is called unduloid,  if $M>0$, $\mathbb{D}(M)$ is not embedded and it is called nodoid. These surfaces are called Delaunays and they are foliated by circles. The trace of the focus in the construction above is called the profile curve of the Delaunay and it is transversal to the circles. Another well known family of constant mean curvature surfaces are the Twizzlers, they are foliated by helices and we can naturally define a profile curve which is transversal to the helices. To make the presentation in this abstract easier, we will be only considering Twizzlers with mean curvature 1. In this paper we will prove that if we roll the profile curve of a Twizzler on a line in a plane and, simultaneously, we move the points in the line at the same speed, so that the rolling motion of the profile curve looks like if it were placed on a treadmill, then, the trace made by the origin of the profile curve is one of the closed heart-shaped curves $\mathbb{H}\mathbb{S}(M,w)=\{(x,y):x^2+y^2+\frac{y}{1+w^2 x^2}=M \}$ for some $M>-\frac{1}{4}$ and $w>0$. We have that, up to a rigid motion in $\bfR{3}$, every Twizzler is determined by its corresponding heart-shaped curve $\mathbb{H}\mathbb{S}(M,w)$; we will denote this Twizzler by $\mathbb{T}(M,w)$. We will prove the continuity of the map $\rho$ from $\Omega=\{M+v\mathfrak{i}\in \mathbb{C}:M\ge -\frac{1}{4},\, M\ne 0,\,  \, 0\le v<1 \}$ to the set of immersions that sends, every point of the form $-\frac{1}{4}+v\mathfrak{i}$ to a cylinder of radius $\frac{1}{2}$ ; a real point $M$ to $\mathbb{D}(M)$ and, for any $M>-\frac{1}{4}$ and $v>0$, the point $M+v\mathfrak{i}$ to  $\mathbb{T}(M,\sqrt{\frac{1-v}{v}})$. This map $\rho$ is one-to-one in the interior of $\Omega$.
Surfaces in the image under $\rho$ of the purely imaginary points in $\Omega$, that is, Twizzlers associated with $M=0$, contain the axis of symmetry and they are special in the sense that each one of these Twizzlers is isometric to a nodoid and to an undoloid, more precisely: the undoloid $\rho(-\frac{\sqrt{v}}{(1+\sqrt{v})^2})$, the twizzler $\rho(v \mathfrak{i})$ and the nodoid $\rho(\frac{\sqrt{v}}{(1-\sqrt{v})^2})$ are isometric. Moreover, we prove that for any $u\in (0,1)$, all the surfaces in the image of the curve $\alpha_u(M)=\rho(M+ \mathfrak{i} \frac{\sqrt{1+4M}-1-2M+u(\sqrt{1+4M}+1+2M)}{\sqrt{1+4M}+1-2M+u(\sqrt{1+4M}-1+2M)}) $ defined in the interval $(-\frac{\sqrt{u}}{(1+\sqrt{u})^2},\frac{\sqrt{u}}{(1-\sqrt{u})^2})$ are isometric. Notice that this curve starts at an undoloid, then it passes through an Twizzler that contains the axis of symmetry and then it ends with at a nodoid. We will prove that every Twizzler that contains the axis of symmetry is contained in the interior of a cylinder of radius 1. Moreover, it is either properly immersed with a discrete group of rotations acting on its isometry group, or it is dense. As an additional result, we prove that if we do the same rolling on a treadmill procedure to the profile curve of a flat surface with the same helicoidal symmetry as the Twizzlers, then, the trace of the origin must lie in a line. We finish the paper by showing the .pdf version generated by the software Mathematica of some of the programs that provide some pictures and animations showing property of Twizzlers and Delanunay surfaces.

\end{abstract}


\maketitle
\section{Introduction}

Let us start this section with some figures, most of them taken from
animations created with Mathematica. The code that produced Figure 1.8 and the videos suggested in Figures 1.3, 1.6, 1.7, and 1.15 will be included at the end of the paper. The reason why this code works is developed throughout the paper.


\begin{figure}[h]\label{profile of undoloids}
\centerline{\includegraphics[width=11.25cm,height=5.79cm]{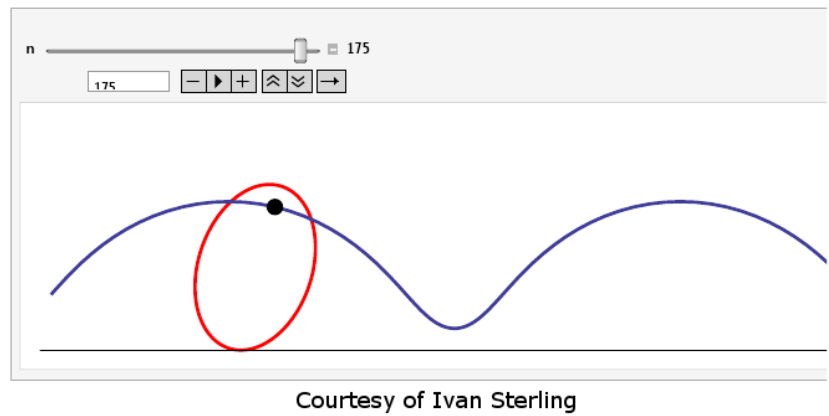}}
\caption{Dynamic interpretation of the profile curve of an unnodoid }
\end{figure}

\vskip2cm

\begin{figure}[h]\label{profile of nodoids}
\centerline{\includegraphics[width=11.25cm,height=5.79cm]{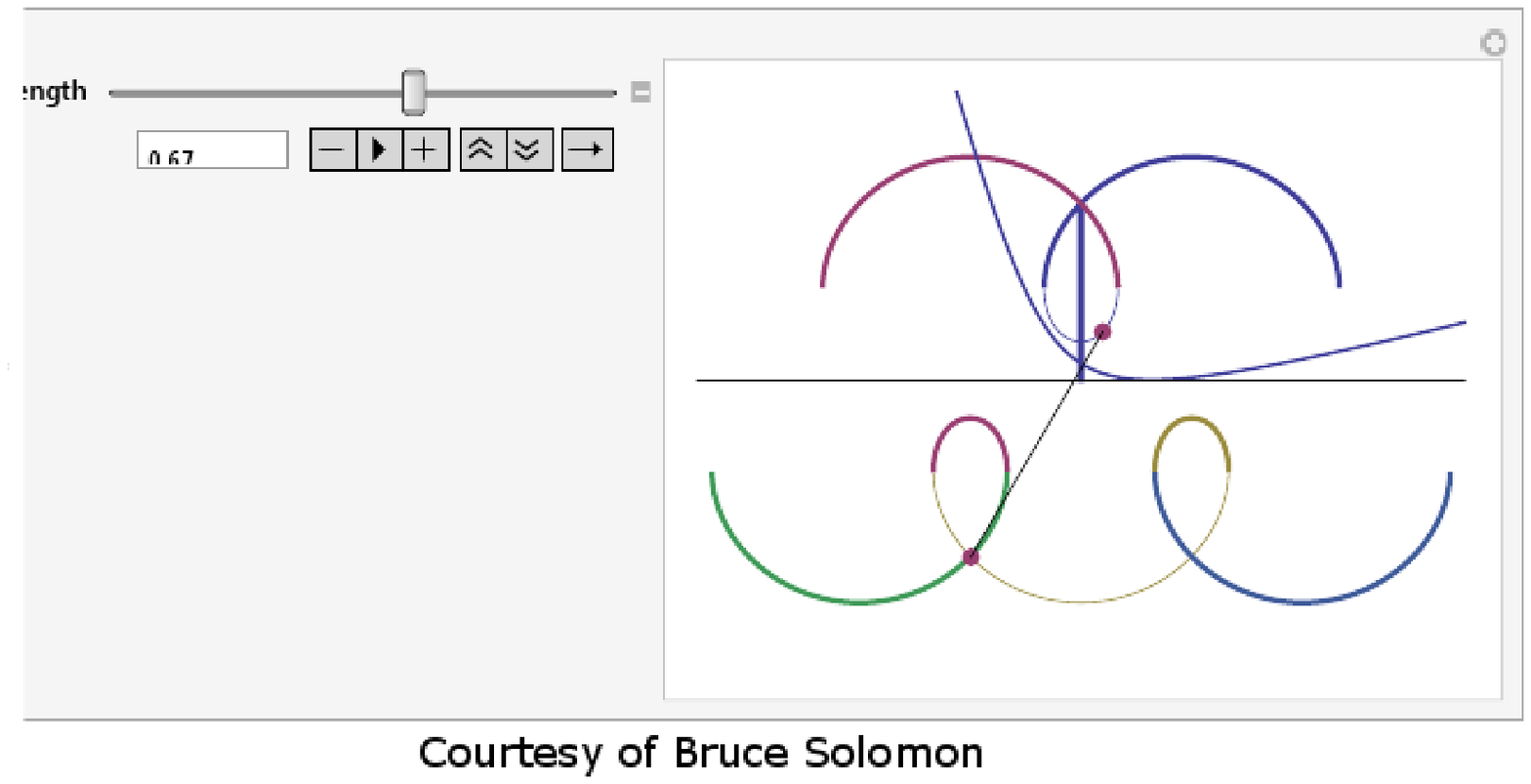}}
\caption{Dynamic interpretation of the profile curve of a nodoid }
\end{figure}

\vfill
\eject

\begin{figure}[h]\label{treadmill ellipse}
\centerline{\includegraphics[width=15cm,height=6.67cm]{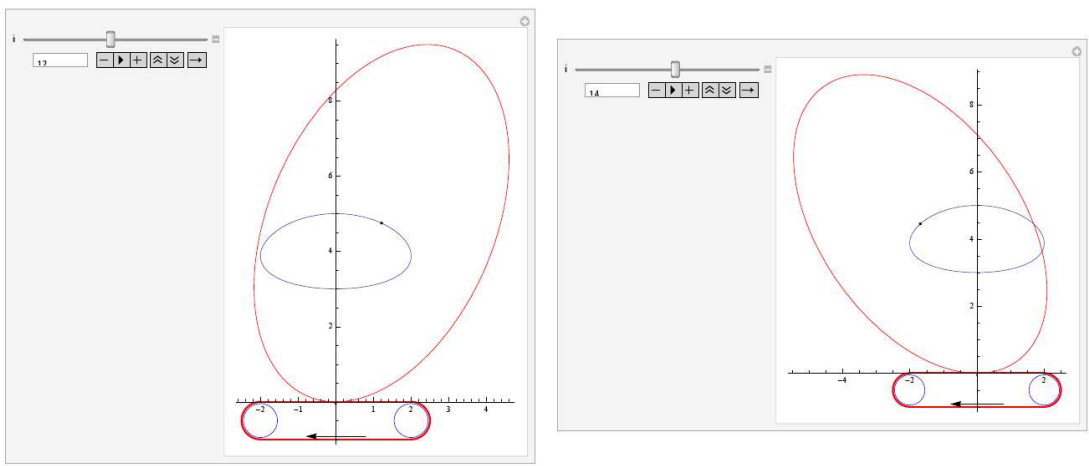}}
\caption{TreadmillSled of an ellipse centered at the origin}
\end{figure}
\vskip.5cm

\begin{figure}[h]\label{contour}
\centerline{\includegraphics[width=15cm,height=4.91cm]{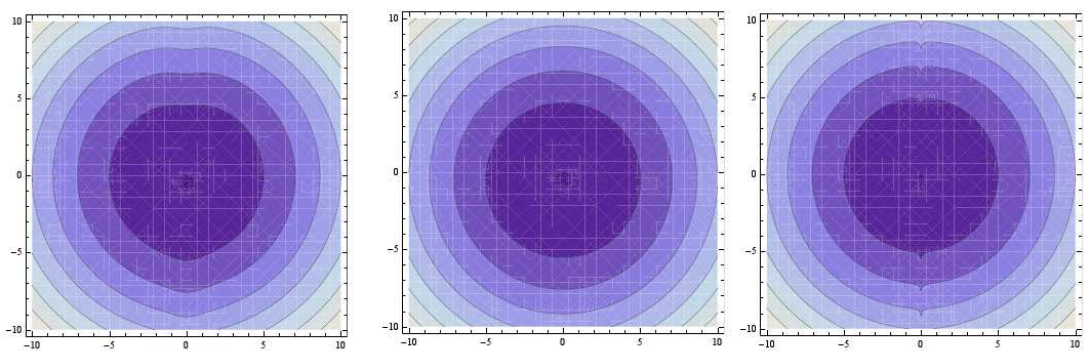}}
\caption{Contours of the integral function $h_w(x,y)=x^2+y^2+\frac{y}{1+w^2x^2}$ for different values of $w$}
\end{figure}
\vskip.5cm

\begin{figure}[h]\label{Profile curve Twizzler}
\centerline{\includegraphics[width=15cm,height=5.31cm]{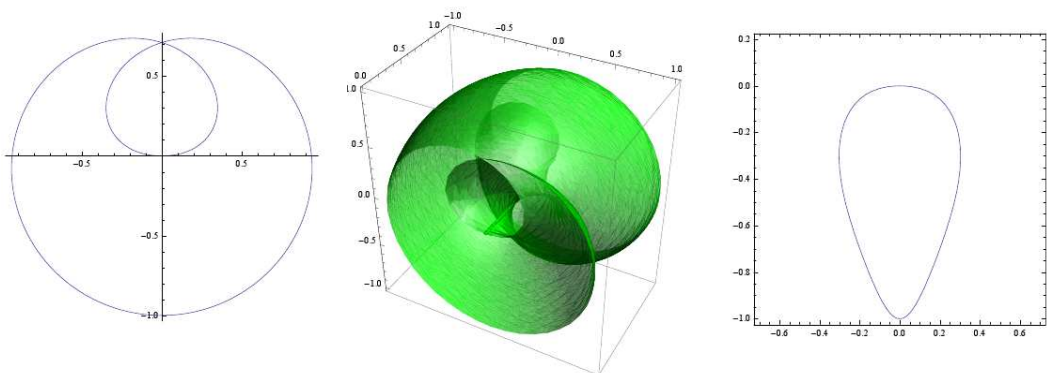}}
\caption{The profile curve, Twizzler and the heart-shaped curve}
\end{figure}
\vfill
\eject

\begin{figure}[h]\label{TreadmillSled of Twizzler}
\centerline{\includegraphics[width=15cm,height=5.00cm]{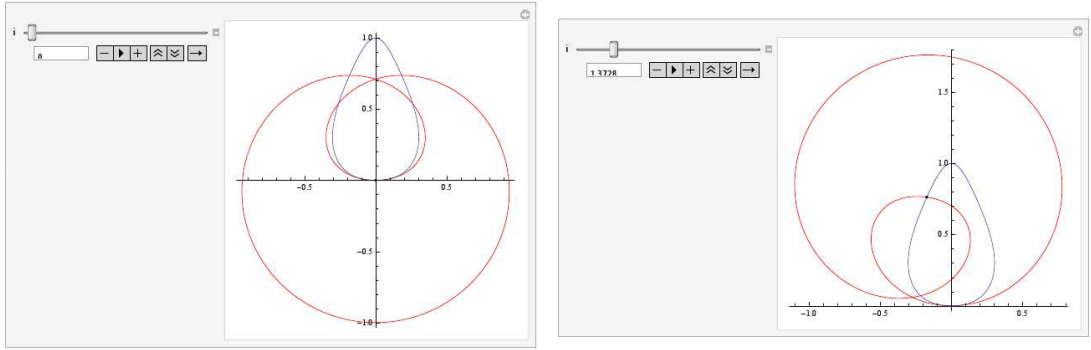}}
\caption{TreadmillSled of the profile curve of a Twizzler}
\end{figure}

\vskip.5cm

\begin{figure}[h]\label{Animation of Moduli space}
\centerline{\includegraphics[width=15cm,height=5.36cm]{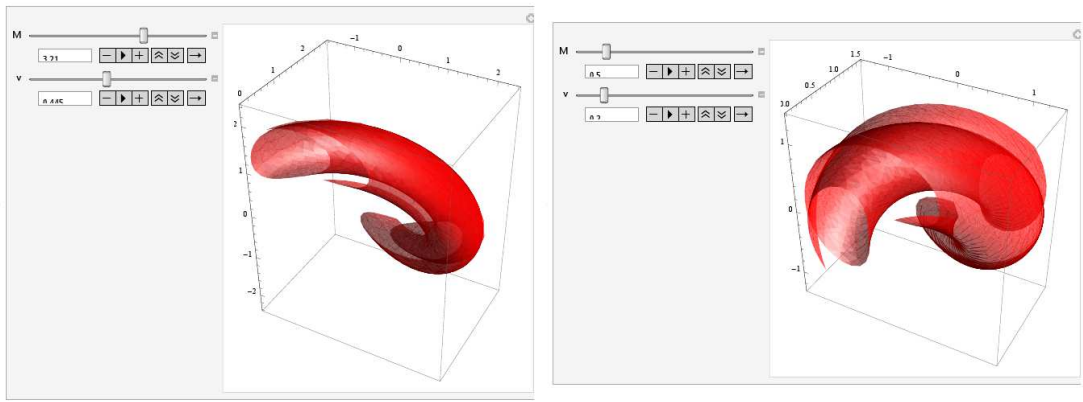}}
\caption{Moduli space of  Twizzler with mean curvature 1.}
\end{figure}

\vskip.5cm

\begin{figure}[h]\label{Fundamental piece}
\centerline{\includegraphics[width=15cm,height=5.45cm]{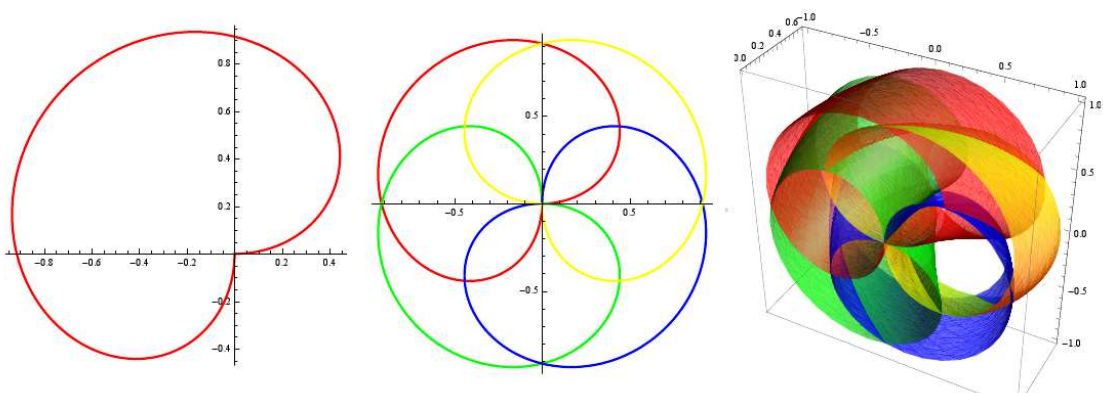}}
\caption{Fundamental piece of the profile curve}
\end{figure}

\vfill
\eject
\begin{figure}[h]\label{profile of surfaces with line}
\centerline{\includegraphics[width=15cm,height=7.14cm]{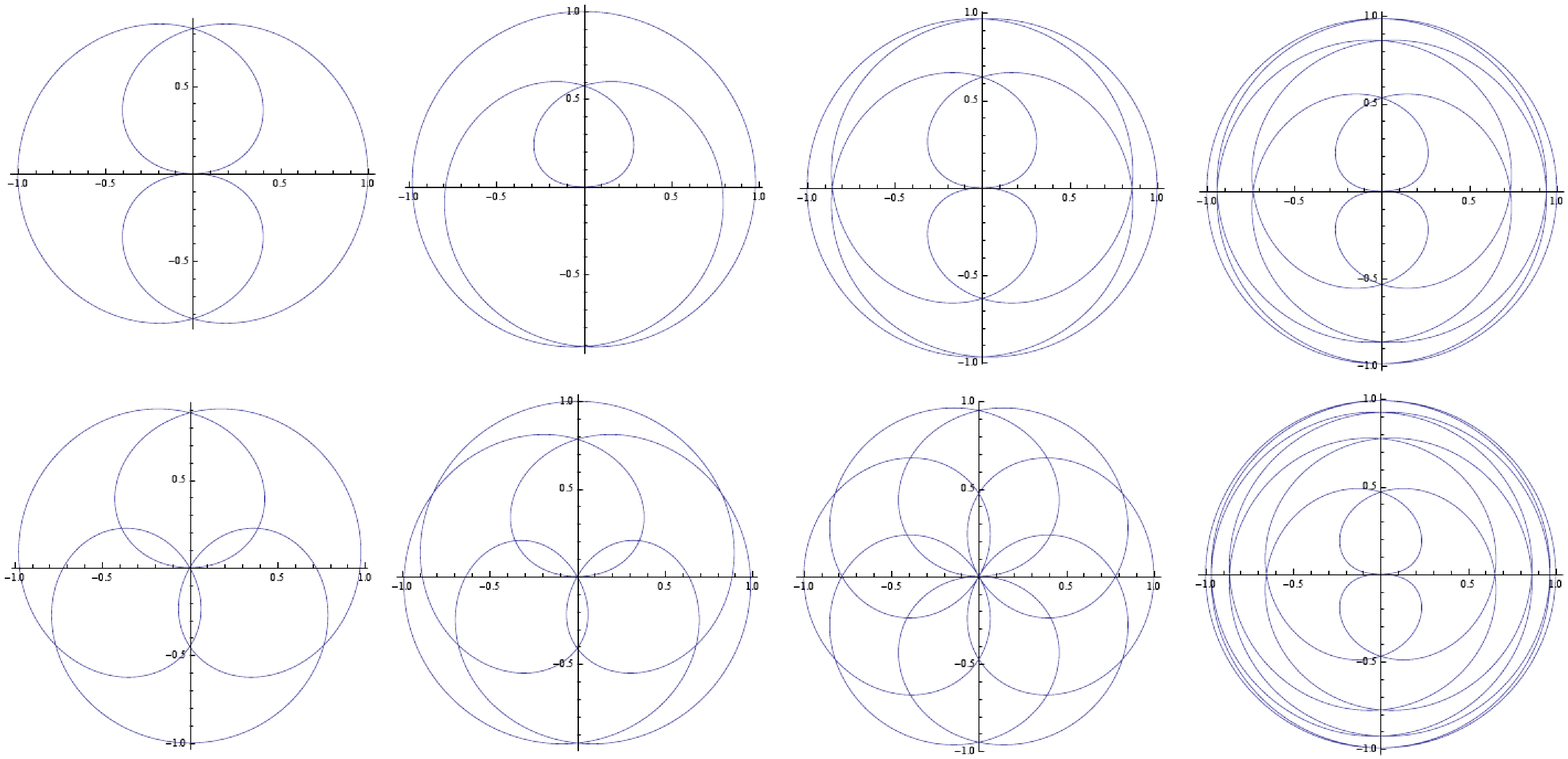}}
\caption{Profile curves of properly immersed Twizzlers that contain the axis}
\end{figure}

\begin{figure}[h]\label{profile of surfaces without line}
\centerline{\includegraphics[width=15cm,height=11.81cm]{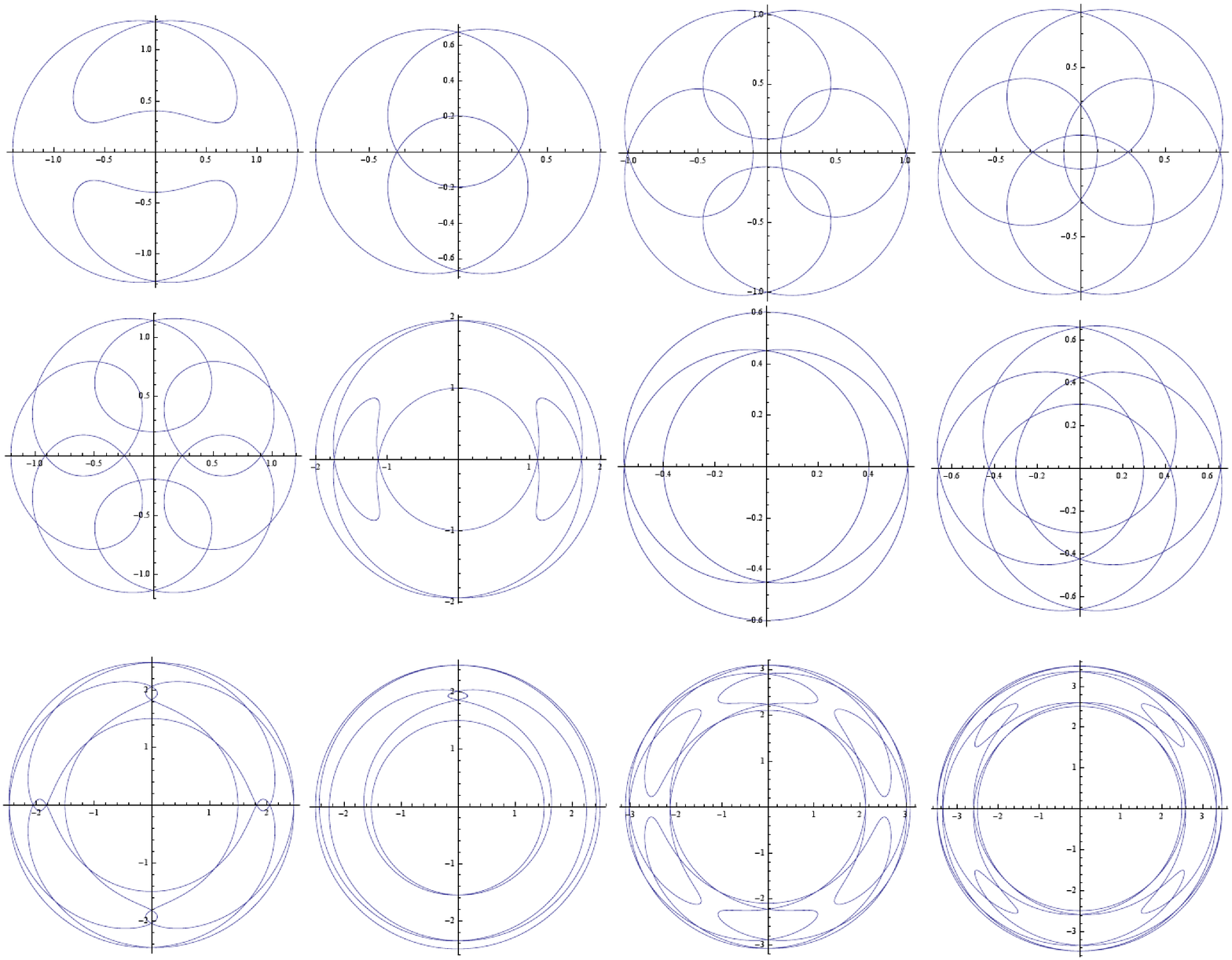}}
\caption{Profile curves of properly immersed Twizzlers that do not contain the axis}
\end{figure}

\vfill
\eject

\begin{figure}[h]\label{profile of dense surfaces}
\centerline{\includegraphics[width=15cm,height=11.81cm]{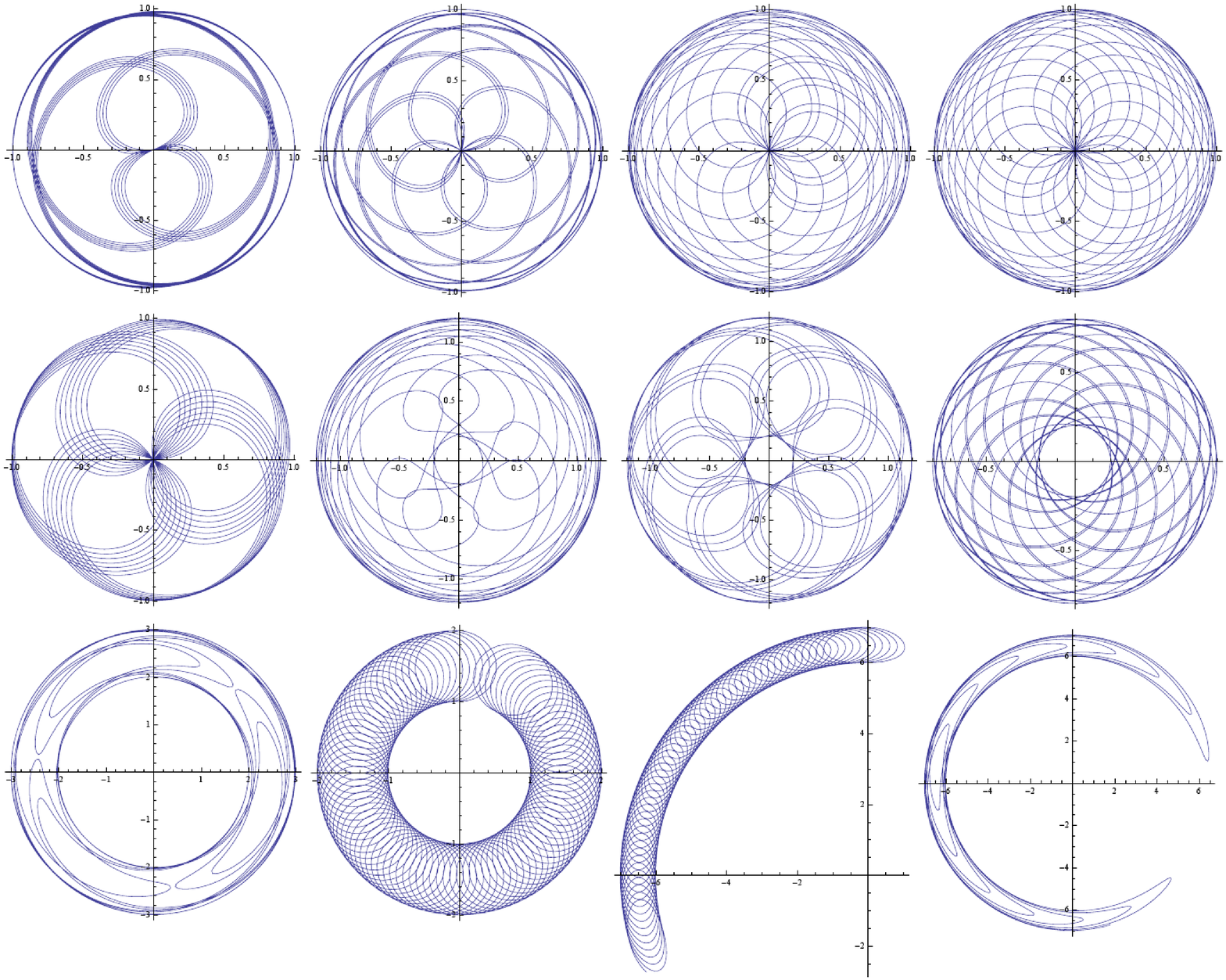}}
\caption{Profile curves of non-properly immersed Twizzlers}
\end{figure}

\vskip1.5cm

\begin{figure}[h]\label{radii}
\centerline{\includegraphics[width=9cm,height=6.5cm]{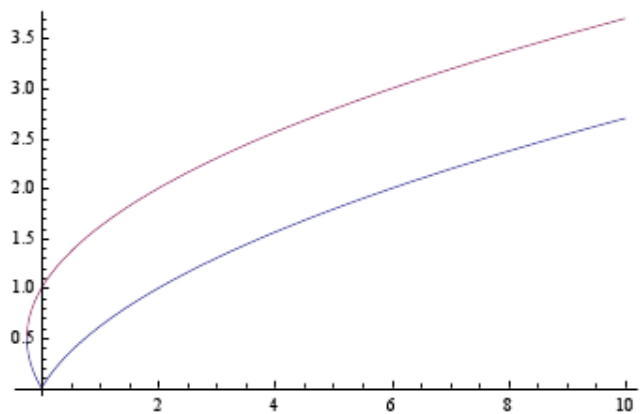}}
\caption{This graph shows how the maximum and minimum distance from the origin to the profile curve of $\mathfrak{T}( M ,v )$ changes with respect to $M$}
\end{figure}

\vfill
\eject

\begin{figure}[h]\label{level set zeta}
\centerline{\includegraphics[width=18cm,height=9cm]{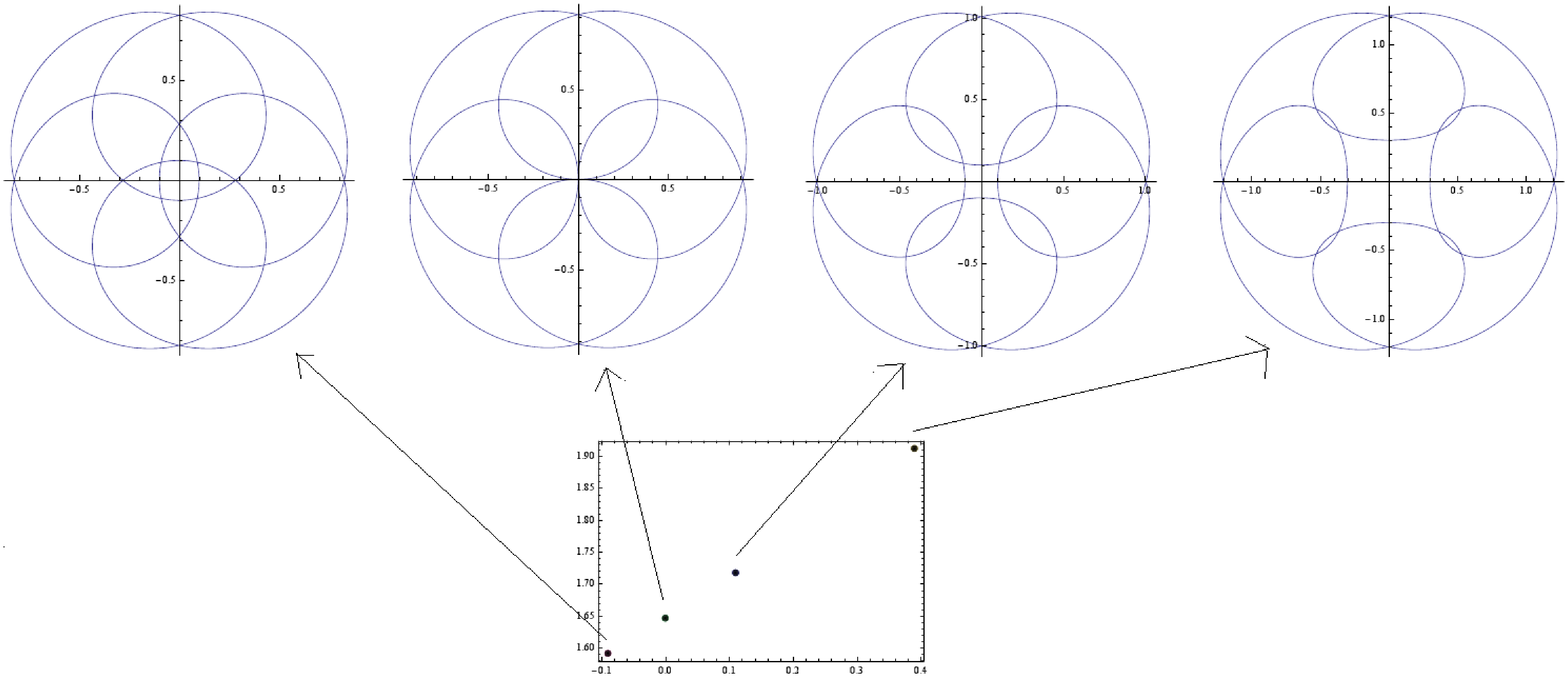}}
\caption{Some profiles curves of Twizzlers that consist of four fundamental pieces and their corresponding values $v$ and $w$}
\end{figure}

\vskip1cm

\begin{figure}[h]\label{Animation of isometric twizzlers}
\centerline{\includegraphics[width=15cm,height=5.56cm]{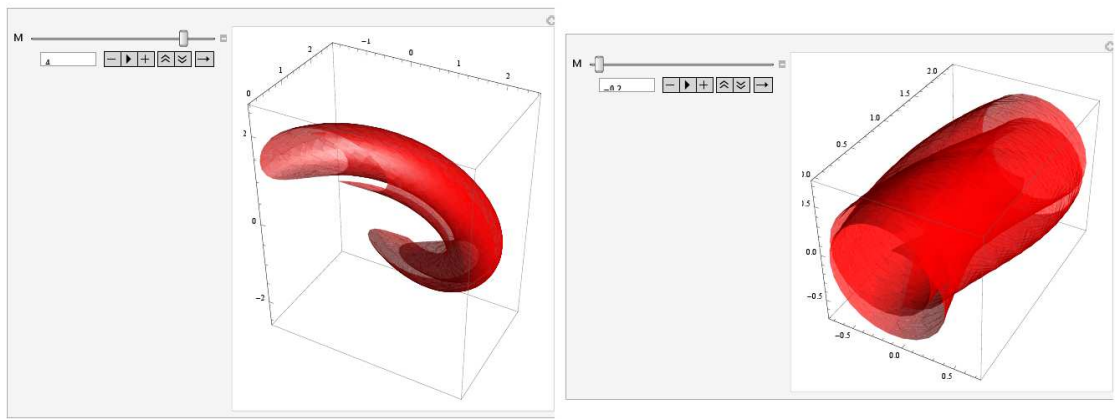}}
\caption{Isometric Twizzlers}
\end{figure}
\vfill
\eject

\vskip1cm

\begin{figure}[h]\label{Animation of some }
\centerline{\includegraphics[width=13cm,height=20.68cm]{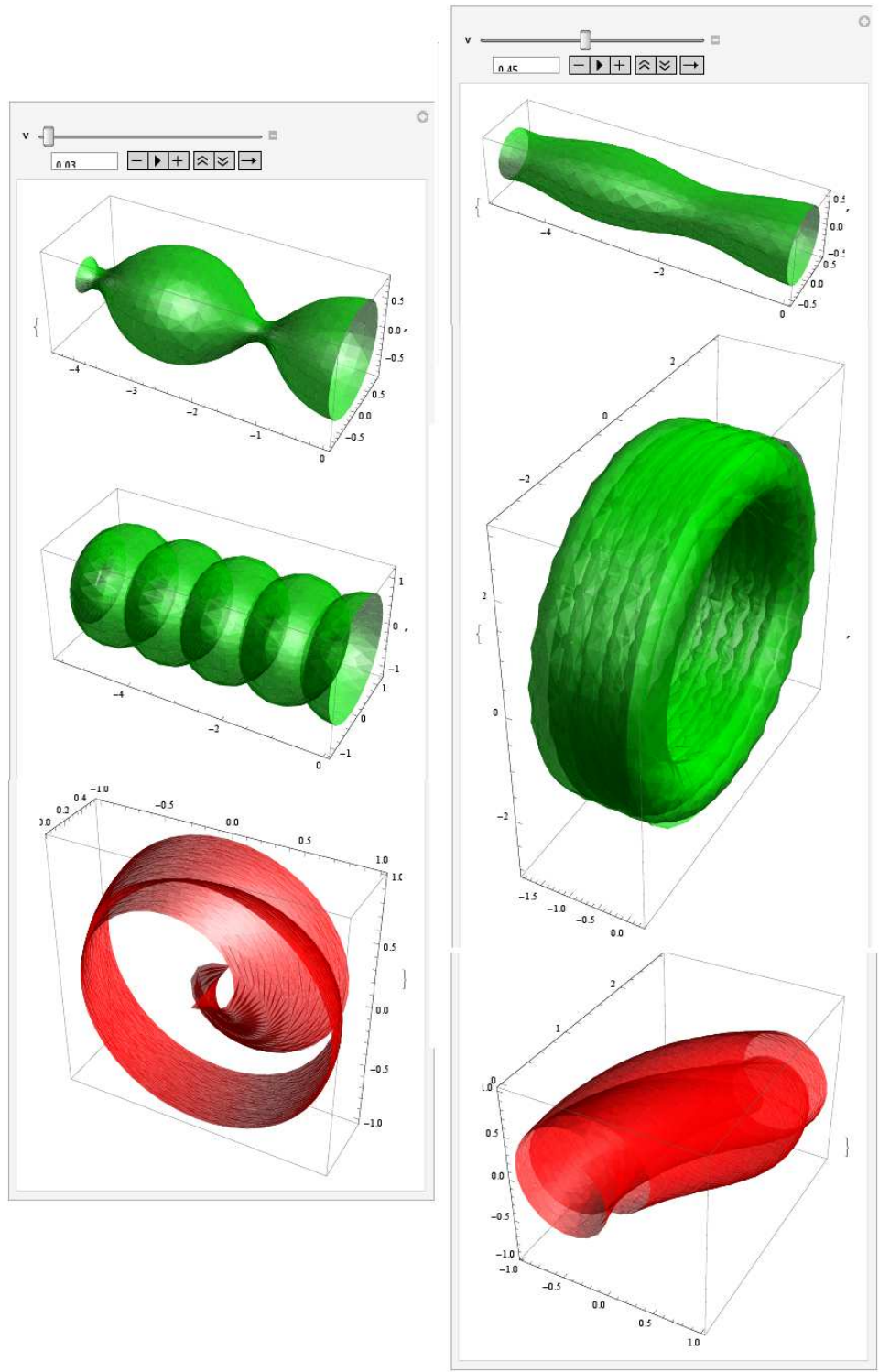}}
\caption{Isometric nodoid, undoloid and special Twizzler}
\end{figure}

\vfill
\eject

\begin{figure}[h]\label{Animation of isometric twizzlers}
\centerline{\includegraphics[width=15cm,height=2.18cm]{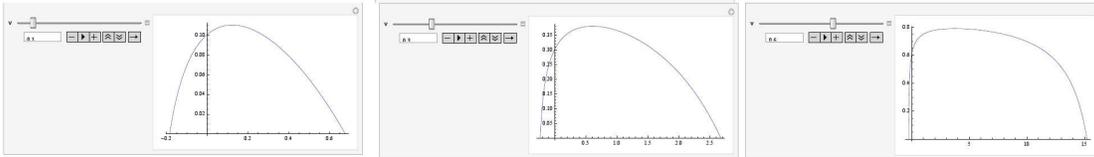}}
\caption{Points in the moduli space that represent isometric Twizzlers}
\end{figure}

\vskip1cm
\begin{figure}[h]\label{radii}
\centerline{\includegraphics[width=15cm,height=5.9cm]{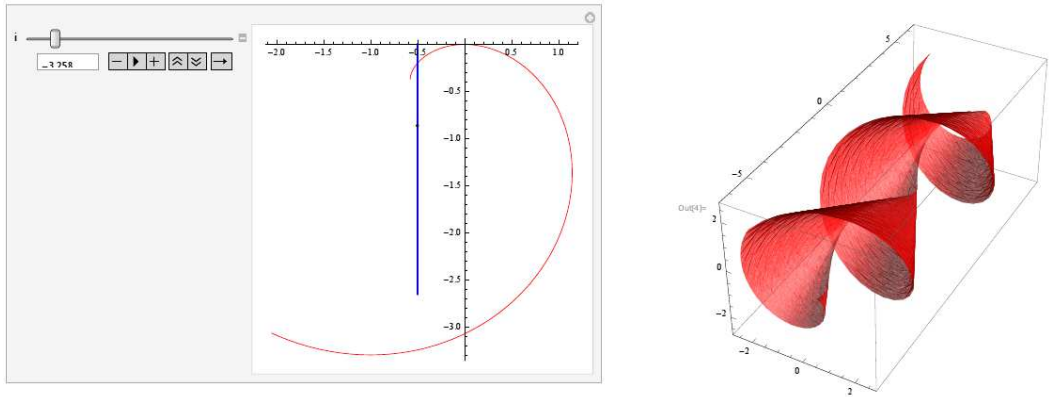}}
\caption{A surface with helicoidal symmetry is flat when its treadmillSled is a vertical semi line}
\end{figure}

In 1841 \cite{De}, Delaunay showed that  if one roll a conic section on a line in a plane and then rotates about that line the trace of a focus, one obtains a cmc surface of revolution. When the conic is a parabola we obtaine a Catenoid, when the conic is an ellipse, the surface is embedded and it is called an {\it  undoloid} and when the conic is a hyperbola the surface is not embedded and it is called a {\it nodoid}. The undoloids and nodoids are called {\it Delaunay surfaces}.  Figure 1.1 illustrates the relation between the ellipse and the trace of its focus. Notice that only one focus is used to get the curve that needs to get rotated in order to generate an undoloid. Figure 1.2 illustrates the relation between the hyperbola and the trace of the foci. Notice that both foci are used to get the curve that needs to get rotated in order to generate a nodoid.

In 1970, Using the integrability of the Gauss Equation and Mainardi-Codazzi equation, Lawson showed that for any immersion $f_0:U\to \bfR{3}$ with constant mean curvature $H$, defined in a simple connected surface $U$, there exists  a
$2\pi$-periodic 1-parametric family of immersions $\{f_{\theta}:U\to \bfR{3}:\theta\in[0,2 \pi]\}$ with constant mean curvature $H$ and with the same induced metric.  The family $f_{\theta}$ is called the {\it associated} family to $f_0$.

\begin{rem}\label{Continuity f theta}
The map $\theta\longrightarrow f_{\theta}$ is continuous with respect to the parameter $\theta$.
\end{rem}

In 1982 \cite{DD}  Do Carmo and Dajczer studied surfaces with constant mean curvature that are invariant under the group $g_t:\bfR{3}\to \bfR{3}$ of rigid motions

$$g_t(x,y,z)=(x\cos(t)+y\sin(t),-x\sin(t)+y\cos(t),z+ht)\quad h \in(-\infty,\infty)$$

They called these surfaces {\it helicoidal} with a pitch $h$. When $h=0$, the group $g_t$ becomes a group of rotations and the helicoidal surfaces become surfaces of revolution. Therefore, when $h=0$ the helicoidal surfaces are either cylinders, spheres or Delaunay surfaces. The main result in their paper established that every Helicoidal surface is associated with a Delaunay surface.

A {\it Twizzler} is an immersion of the form

\begin{eqnarray}\label{surfaces}
\phi(s,t)=(x(s) \cos(w t)+z(s) \sin(w t),\, t\, ,-x(s) \sin(w t)+z(s) \cos(w t))
\end{eqnarray}
 with constant mean curvature. We will assume that the curve $(x(s), z(s))$ is parametrized by arc length and we will call it the {\it profile curve} of the Twizzler. Notice that Twizzlers correspond to those helicoidal surfaces with nonzero pitch.

In this paper we will give an interpretation of the profile curve of Twizzlers similar to the interpretation for the profile curves of the Delaunays. In order to do this,
given a regular curve $\alpha$ in the plane, we define the {\it TreadmillSled} of $\alpha$ as the trace of the origin when this curve (and the plane that contains it) rolls on a treadmill located at the origin which is aligned along the $x$-axis. Notice that the origin, as a point in the plane that contains the curve $\alpha$, plays and important roll in this definition, for example, the TreadmillSled of a circle with center at the origin is just a point, while the TreadmillSled of a circle of radius $R$ whose center is at a distance $r$ of the origin, is a circle of radius $r$ with center $(0,R)$. Figure 1.3 shows the TreadmillSled of an ellipse with center at the origin. The black dot represents the center of the ellipse while the blue curve represents the trace of this center. As the next theorem shows, it turns out that for any curve $(x(s),z(s))$ parametrized by arc-length there is an easy expression for its TreadmillSled.

\begin{thm}\label{Expression for TreadmillSled}
If  $\alpha(s)=(x(s),z(s))$ is a curve parametrized by arc-length and

$$\xi_1(s)=x(s) \, x^\prime(s)+z(s)\, z^\prime(s)\com{and} \xi_2(s)=-x(s)\, z^\prime(s)+z(s)\, x^\prime(s)$$

then the TreadmillSled of the curve $\alpha$ is $-(\xi_1(s),\xi_2(s))$
\end{thm}

Notice that $(\xi_1(s),\xi_2(s))$ are the coordinates of the vector $(x(s),z(s))$ with respect to the orthonormal basis $\{(x^\prime(s),z^\prime(s)),(-z^\prime(s),x^\prime(s))\}$.

There is a relation between the profile curve of Delaunay surfaces and conics. For Twizzlers, there is a relation between their profile curves and the level sets of the function $h_w(x,y)=x^2+y^2+\frac{y}{\sqrt{1+w^2 x^2}}$. It is not difficult to check that the range of the function $h_w$ is the interval $[-\frac{1}{4},\infty)$, that $h_w^{-1}(-\frac{1}{4})=\{(0,-\frac{1}{2})\}$ and that, every $M>-\frac{1}{4}$ is a regular value of $h_w$ and $h_w^{-1}(M)$ is a closed simple curve. We will refer to these level sets as {\it heart-shaped} curves. Figure 1.4 shows some of these level sets

The following Theorem establish the relation between the heart-shaped curves and the profile curve of the Twizzlers

\begin{thm}\label{Profile Twizzler vs heart-shape curves}
The treadmillSled of the profile curve of a Twizzler with constant mean curvature 1 is a heart-shaped curve $-h_w^{-1}(M)$ for some $M\ge-\frac{1}{4}$. The value $M=-\frac{1}{4}$ is achieved by a cylinder of radius $\frac{1}{2}$.
\end{thm}

The previous theorem is a consequence of the following proposition:

\begin{prop}\label{ODE for the natural coordinates} The immersions given by (\ref{surfaces}) have mean curvature 1 if and only if the functions $\xi_1$ and $\xi_2$ defined in Lemma \ref{Expression for TreadmillSled} satisfy the following ordinary differential equations $\xi_1^\prime(s)=f_1(\xi_1(s),\xi_2(s))$ and $\xi_2^\prime(s)=f_2(\xi_1(s),\xi_2(s))$, where,

\begin{eqnarray}
f_1(x_1,x_2)&=&\frac{-w^2 x_2 + 2(1+w^2 x_1^2)^\frac{3}{2}}{1+w^2(x_1^2+x_2^2)} \, x_2+1\cr
& &\cr
f_2(x_1,x_2)&=&\frac{w^2 x_2 - 2(1+w^2 x_1^2)^\frac{3}{2}}{1+w^2(x_1^2+x_2^2)}\, x_1
\end{eqnarray}

Moreover, the function $h_w$ is an integral curve for every solution $(\xi_1(s),\xi_2(s))$.

\end{prop}

Figure 1.5 shows the profile curve of a Twizzler and the heart-shaped curve associated with it. Figure 1.6, is a picture showing how, for this Twizzler,  the TreadmillSled of the profile curve is indeed negative the heart-shaped curve.

The following theorem gives us another application of the treadmillSled

\begin{thm}\label{Flat surfaces} A surfaces of the form (\ref{surfaces}) is flat if and only if either the TreadmillSled of the profile curve is a point in the $y$-axis different from the origin, in this case the surface is a cylinder, or the treadmillSled of the profile curve is contained in a vertical semi line that starts at point in the $x$-axis different from the origin. Moreover, the functions  $x$ and $z$ can be explicitly computed and they are:

$$ x(s)= \frac{1}{2}\cos{\frac{2\sqrt{as+b}}{a}}+\sqrt{as+b}\sin{\frac{2\sqrt{as+b}}{a}}\com{and}
z(s)= \sqrt{as+b}\cos{\frac{2\sqrt{as+b}}{a}} - \frac{1}{2}\sin{\frac{2\sqrt{as+b}}{a}} $$
\end{thm}

Figure 1.17 shows how the TreadmillSled of the profile curve produces a vertical semi line.

If we exclude the cylinder and the value $M=-\frac{1}{4}$, Theorem \ref{Profile Twizzler vs heart-shape curves} establishes a 1:1 correspondence between pairs $(M,w)$ with $M>-\frac{1}{4}$ and $w > 0$ and Twizzlers with mean curvature 1. In order to deal with a bounded parameter, we will define $v = 1/(1+w^2)$, therefore, the parameter $v$ moves from $0$ to $1$ when $w$ moves from $\infty$ to $0$. Figure 1.7 shows pictures from an animation that produces a piece of the Twizzler associated with values of $M$ and $v$. We will refer to this Twizzler as  $\mathfrak{T}( M ,v )$ when the dependence of $M$ and $v$ is needed.

Explicit parametrizations for the heart-shaped curves make the understanding and drawing of Twizzlers easier. The following Lemma provides a possible parametrization for these heart-shaped curves.

\begin{lem}\label{Parametrization heart-shaped} For any $M>-\frac{1}{4}$ and $w>0$, the curve $\alpha(t)=(\rho_1(t),\, \rho_2(t))$ defined on the interval $[0,2\pi]$ and given by

$$\rho_1(u)=A \cos(u)\com{and}
\rho_2(u)=\frac{-1+\sqrt{1+4M+B\cos^2(u)}\, \sin(u)}{2\sqrt{1+w^2A^2\cos^2(u)}}$$

where,

$$A= \frac{\sqrt{-1+Mw^2+\sqrt{1+(1+2M)w^2+M^2w^4}}}{\sqrt{2}\, w}\com{and}$$

$$B= \frac{2+2M^2w^4+w^2+2(Mw^2-1)\sqrt{1+(1+2M)w^2+M^2w^4}}{w^2}$$

is a closed simple regular curve that parametrizes the heart-shaped curve $h_w^{-1}(M)$.
\end{lem}

We can define a {\it fundamental piece of a profile curve} as a connected piece of profile curve with the property that the TreadmillSled motion of this piece goes exactly once over the heart-shaped curve. It is not difficult to see that the whole profile curve is the union of fundamental pieces. Figure 1.8 shows the fundamental piece of the profile of a properly immersed Twizzler, along with the whole profile curve made up of four pieces in this case, and the graph of the Twizzler. For an undoloid we could define a fundamental piece as the trace of the focus of the ellipse when this ellipse rolls once. It is clear that the integers $\mathbb{Z}$ acts on the group of symmetries of the Delaunay surfaces in the form of translations. Theorem \ref{Moduli space Twizzler} shows that the group $\mathbb{Z}$ also acts on the set of isometries of Twizzlers.

In order to produced Figure 1.8, it was necessary to find a formula for the length of the fundamental piece of the profile curve. The following lemma provides such a formula.

\begin{lem}\label{Length fundamental piece}
The length of the fundamental piece of the Twizzler $\mathfrak{T}( M ,v )$ is $\int_0^{2 \pi}\sqrt{\frac{\lambda}{\mu}}du$
where,
\begin{eqnarray*}
\lambda(u)=&
 (\frac{d\rho_1}{du})^2+(\frac{d\rho_2}{du})^2\cr
 & &\cr
 \mu(u)=&f_1^2(\rho_1(u),\rho_2(u))+f_2^2(\rho_1(u),\rho_2(u)).\cr
\end{eqnarray*}
The functions $\rho_1$, $\rho_2$, $f_1$ and $f_2$ are defined in the Lemma \ref{Parametrization heart-shaped} and the Proposition \ref{ODE for the natural coordinates}. Recall that $w$ and $v$ are related by the equation $v = 1/(1+w^2)$.
\end{lem}

The following theorem along with Theorem \ref{Twizzler's properties} and Theorem \ref{Twizzler's moduli space} give a precise picture of the Moduli space for the Twizzlers.

\begin{thm}\label{Moduli space Twizzler} Every Twizzler $\mathfrak{T}( M ,v )$ is invariant under a group of rotations of the form $\{ R(n \theta) \,:\,  n \in \mathbb{Z} \}$, where the angle $\theta$ depends on the fundamental piece of the profile curve. Moreover, if $R(m \theta)=R(\theta)$ for some integer $m$, then the Twizzler is properly immersed, otherwise it is dense in the interior of a cylinder of radius 1 when $M=0$ or it is dense in the region bounded by two concentric cylinders of radii
$r_1(M) =| \frac{\sqrt{1+ 4 M}-1}{2}|$  and $r_2(M) = \frac{\sqrt{1+ 4 M}+ 1}{2}$ in the case that $M\ne 0$. We also have another type of density: the set of points $(M,v)$ associated with properly immersed Twizzler is non countable and dense. In addition, we have that $\mathfrak{T}( M ,v )$ is a properly immersed surface with profile curve consisting of $b$ fundamental pieces, if and only if $\int_0^{2\pi}\psi du=2\pi \frac{a}{b}$ with $a$ and $b$ positive relatively prime integers, where

$$\psi(u)=\frac{-w^2\rho_2(u) + 2(1+w^2\rho_1^2(u))^\frac{3}{2}}{1+w^2(\rho_1^2(u)+\rho_2^2(u))}\, \sqrt{\frac{\lambda(u)}{\mu(u)}}$$

\end{thm}

For a Twizzler that does not contain the axis of symmetry, the ``properly immersed  vs. dense''  property established in the theorem above was proven by Hitt and Roussoss \cite{HR}.  By numerically solving the equation $\int_0^{2\pi}\psi du=2\pi \, a/b$ in  the previous theorem, we can graph profile curves of Twizzlers with any desire property. In Figures 1.9 we solve the numerical equation fixing $M=0$ and taking several integer values for $a$ and $b$. Since $M=0$, these profile curves represent twizzlers that contain the axis of symmetry. In Figure 1.10 we took several values for $M\ne 0$ and $a$ and $b$ integers to produce properly immersed twizzlers that do not contain the axis of symmetry. In Figure 1.11 we took $a$ and $b$ so that $a/b$ is not rational so that the twizzler is not properly immersed.  In Figure 1.13 we took $a=5$, $b=4$ and 4 values of $M$ in order to produce properly immersed examples, we also show the points $(M,w)$ associated with these Twizzlers. Figure 1.12 shows the graph of the function $r_1(M)$ or $r_2(M)$ defined in Theorem \ref{Moduli space Twizzler}.

By Dajczer and Do Carmo's Theorem, we know that every Twizzler is isometric to a Delaunay surfaces, in order to decide which one, we will prove the following properties for  Delaunay surfaces,

\begin{thm}\label{Delaunay's properties}
For every nonzero real number  $M\in  (-\frac{1}{4},\infty)$, the Delaunay surface  $\mathbb{D}(M)$ generated by the conic $\{(x,y): 4 x^2-\frac{y^2}{M}=1\}$ has constant mean curvature 1. Moreover, the quotient between the maximum value of the Gauss curvature  and the minimum value of the Gauss curvature of  $\mathbb{D}(M)$ is given by the function $rs(M)=-(\frac{1-\sqrt{1+4 M}}{1+\sqrt{1+4 M}})^2$

\end{thm}

We have that the function $rs$ defines a bijection between the intervals $(-\frac{1}{4},0)$ to $(0,1)$ and also it  defines a bijection between the intervals $(0,\infty)$ to $(0,1)$. On the other hand, each undoloid is isometric (associated with angle $\theta=\frac{\pi}{2}$) to a nodoid, see for example Appendix A in the paper \cite{K} written by  Kopouleas. As a consequence
of Theorem \ref{Delaunay's properties} we have the following corollary:

\begin{cor} \label{rate max min Gauss curvature}
Two helicoidal surfaces are isometric, if and only if the quotient of the maximum value and the minimum value of the Gauss curvature is the same, in particular, for any $u\in (0,1)$ the undoloid $\mathbb{D}(-\frac{\sqrt{u}}{(1 + \sqrt{u})^2})$ is isometric to the nodoid $\mathbb{D}(\frac{\sqrt{u}}{(1 - \sqrt{u})^2})$
\end{cor}

As pointed out before, each nodoid is isometric to a undoloid, therefore we can replace the word Delaunay by either the word undoloid or nodoid in the  Do Carmo-Dajczer theorem, that is, we can say that each Twizzler is isometric (or associated) to either a nodoid or an undoloid. Another family of surfaces that holds the same property is the set of Twizzlers that contains the axis of symmetry, that is, the set of Twizzlers corresponding to $M=0$ in the moduli space. We will call these surfaces {\it special twizzlers} and we will denote them by $\mathbb{S}\mathfrak{T}(v)$, that is, $\mathbb{S}\mathfrak{T}(v)=\mathfrak{T}( 0 ,v )$.

The following Theorem gives us the quotient of the maximum value and the minimum value of the Gauss curvature for special Twizzlers.

 \begin{thm}\label{special Twizzler's properties}
For every nonzero real number  $v\in  (0,1)$, the quotient between the maximum value of the Gauss curvature  and the minimum value of the Gauss curvature of the special Twizzler surface  $\mathbb{S}\mathfrak{T}(v)$, is $-v$. Moreover $\mathbb{S}\mathfrak{T}(v)$ is isometric to the undoloid $\mathbb{D}(-\frac{\sqrt{v}}{(1 + \sqrt{v})^2})$ and the nodoid $\mathbb{D}(\frac{\sqrt{v}}{(1 - \sqrt{v})^2})$.
\end{thm}

Figure 1.15 shows two sets of isometric nodoid-undoloid-special Twizzler surfaces. We can generalize the theorem above as follows,

\begin{thm}\label{Twizzler's properties}
If $v=\frac{1}{1+w^2}$, the quotient between the maximum value of the Gauss curvature  and the minimum value of the Gauss curvature of the Twizzler surface $\mathfrak{T}( M,v)$ is $-\frac{2+(1+2M-\sqrt{1+4M})\, w^2}{2+(1+2M+\sqrt{1+4M})\, w^2}$. Moreover, fixing $c\in(0,1)$, all the Twizzlers in the set

$$\{\mathfrak{T}( M,\frac{\sqrt{1+4M}-1-2M+c(\sqrt{1+4M}+1+2M)}{\sqrt{1+4M}+1-2M+c(\sqrt{1+4M}-1+2M)}) \, :\,
M\in (-\frac{\sqrt{c}}{(1+\sqrt{c})^2},\frac{\sqrt{c}}{(1-\sqrt{c})^2}) \}$$

are isometric.
\end{thm}

Since the curve $\alpha_c(M)=(M, \frac{\sqrt{1+4M}-1-2M+c(\sqrt{1+4M}+1+2M)}{\sqrt{1+4M}+1-2M+c(\sqrt{1+4M}-1+2M)})$ satisfies that
$\alpha_c(-\frac{\sqrt{c}}{(1+\sqrt{c})^2})=(-\frac{\sqrt{c}}{(1+\sqrt{c})^2},0)$, $\alpha_c(0)=(0,c)$ and $\alpha_c(\frac{\sqrt{c}}{(1-\sqrt{c})^2})=(\frac{\sqrt{c}}{(1-\sqrt{c})^2},0)$\, then,
we have as a corollary of the previous two theorems and Remark \ref{Continuity f theta} that

\begin{thm}\label{Twizzler's moduli space}
Let $\Omega=\{M+\mathfrak{i} v\in\mathbb{C}: M\ge-\frac{1}{4},\, M\ne0\com{and} 0\le v<1\}$. The function $\rho$ from $\Omega$ to the set of immersions in $\bfR{3}$ given by

\begin{eqnarray*}
\rho(M+\mathfrak{i} v)&=&\mathfrak{T}( M,v) \com{for any $v>0$ and $M\ne -\frac{1}{4}$}\cr
\rho(M)&=&\mathbb{D}(M) \com{for any $M\ne 0$ and $M\ne -\frac{1}{4}$}\cr
\rho(-\frac{1}{4}+\mathfrak{i} v)&=&\{ (x,y,z)\in\bfR{3}: x^2+z^2=\frac{1}{4} \}\cr
\end{eqnarray*}

is continuous in the sense that for every point $p$ in $\Omega$ there exists a neighborhood $U$ of $p$ in $\Omega$ and a continuous function $f:U\times \bfR{2}\longrightarrow \bfR{3}$
such that for any
$M+\mathfrak{i} v\in U$, the map  $(s,t)\longrightarrow f(M+\mathfrak{i} v,s,t)$ defines a parametrization of the surface
$\rho(M+\mathfrak{i} v)$. Moreover, the function $\rho$ is one to one in the interior of $\Omega$.
\end{thm}

The continuity at the points of form $-\frac{1}{4}+\mathfrak{i} v$ follows from Theorem \ref{Moduli space Twizzler}, because each Twizzler $\mathfrak{T}( M,v)$ is contained the region bounded by the two concentric cylinders of radii
$r_1(M) =| \frac{\sqrt{1+ 4 M}-1}{2}|$  and $r_2(M) = \frac{\sqrt{1+ 4 M}+ 1}{2}$. Figure (1.17) shows the trace of the curve $\alpha_u$ for several values of $u$. As a consequence of the previous theorem we have that the family of associated surfaces $f_{\theta}$ for Twizzlers repeats every $\frac{\pi}{2}$. More precisely, it is $\pi$-periodic and it is 2:1 when defined in the interval $[0,\pi]$.

The author would like to express his gratitude to Professors Robert Kusner, Ivan Sterling, Bruce Solomon, Wayne Rossman,
Martin Kilian and Ioannis Roussos for solving several of his doubts about Delaunays and Twizzlers and for providing him references.

\section{Proof of the results}

In this section we will provide the proofs of the results established in the introduction. We will start with Theorem \ref{Expression for TreadmillSled}.

{\bf Proof Theorem \ref{Expression for TreadmillSled}:}  For every $s_0$ in the domain of the curve $\alpha$, we can find an isometry from $\bfR{2}$ to  $\bfR{2}$ of the form $T{\bf x}= R(\gamma(s_0)){\bf x}+c(s_0)$ that rotates the curve $\alpha$ so that $\alpha(s_0)$ moves to the origin and the vector $\alpha^\prime(s_0)$ moves to the unit vector that points to the positive direction of the $x$ axis. Here

$$R(t)=\begin{pmatrix}\cos(t) & -\sin(t)\cr \sin(t)& \cos(t)\end{pmatrix}$$

The angle $\gamma(s_0)$ and the constant $c(s_0)$ can be found by the solving the following equations (the solution of these equations is the key part of the proof of this theorem)

$$\beta^\prime(s_0)=\begin{pmatrix}1\cr 0\end{pmatrix}\com{and } \beta(s_0)=\begin{pmatrix}0\cr 0\end{pmatrix}$$

where $\beta(s)=R(\gamma(s_0))\alpha(s)+c(s_0)$. If follows that the TreadmillSled of the curve $\alpha$ is the curve $c(s)$. A direct verification shows that

$$c(s)=\begin{pmatrix} -x^\prime(s)x(s) -z^\prime(s)z(s)\cr
z^\prime(s)x(s)-x^\prime(s)z(s)\end{pmatrix}$$

This completes the proof of the theorem. \eop

{\bf Proof of Proposition \ref{ODE for the natural coordinates}:} Since the curve $(x(s),z(s))$ is parametrized by arc-length, we can consider
a function  $\theta(s)$ such that

$$x^\prime(s)=\cos(\theta(s))\quad\hbox{and} \quad z^\prime(s)=\sin(\theta(s))$$

Let us define the functions $\xi_1(s)$ and $\xi_2(s)$ as

$$\xi_1=x \cos(\theta)+z\sin(\theta)\quad \hbox{and} \quad\xi_2=-x \sin(\theta)+z\cos(\theta) $$

A direct verification shows that

\begin{eqnarray}\label{x and z in terms of xi1 and xi2}
x=\xi_1 \cos(\theta)-\xi_2\sin(\theta),\quad \quad z=\xi_1 \sin(\theta)+\xi_2\cos(\theta)\quad\hbox{and} \quad \theta^\prime=x^\prime z^{\prime\prime}-z^\prime x^{\prime\prime}
\end{eqnarray}

Moreover, it is not difficult to check that,

$$\xi_1^\prime=\theta^\prime \xi_2+1,\quad \xi_2^\prime=-\theta^\prime \xi_1 \com{and} \xi_1^2+\xi_2^2=x^2+z^2$$

A direct verification shows that the first fundamental form of $\phi$ is given by

$$ E=\<\phi_s,\phi_s\>=1,\, F=\<\phi_s,\phi_t\>=w(zx^\prime-xz^\prime)=w\xi_2, \, G=\<\phi_t,\phi_t\>=1+w^2(x^2+z^2)=1+w^2(\xi_1^2+\xi_2^2)$$

and therefore,

$$EG-F^2=1+w^2(\xi_1^2+\xi_2^2)-w^2\xi_2^2=1+w^2\xi_1^2\, .$$

A Gauss map of the immersion $\phi$ is given by $\nu=\frac{1}{\sqrt{EG-F^2}}\, \phi_s\times \phi_t$. A direct verification shows that,

$$\nu(s,t)=\frac{1}{\sqrt{1+w^2\xi_1^2(s)}}\, (\sin(w t-\theta(s)),w\xi_1 ,\cos(wt-\theta(s)))$$

Again, a direct verification shows that the second fundamental form of $\phi$ is given by

$$ e=\<\phi_{ss},\nu\>= \frac{\theta^\prime}{\sqrt{1+w^2\xi_1^2}}  ,\quad f=\<\phi_{st},\nu\>=\frac{-w}{\sqrt{1+w^2\xi_1^2}},\com{and} g=\<\phi_{tt},\nu\>=\frac{-w^2\xi_2}{\sqrt{1+w^2\xi_1^2}}  $$

Therefore, if we assume that the mean curvature $\frac{eG-2fF+gF}{2(EG-F^2)}$ is 1 we obtain the following ODE

\begin{eqnarray}\label{derivative theta}
\theta^\prime=\frac{-w^2\xi_2 + 2(1+w^2\xi_1^2)^\frac{3}{2}}{1+w^2(\xi_1^2+\xi_2^2)}
\end{eqnarray}

Using the expression above for $\theta^\prime$ in the equations $\xi_1^\prime=\theta^\prime \xi_2+1,\quad \xi_2^\prime=-\theta^\prime \xi_1$, we obtain that $\xi_1$ and $\xi_2$ satisfy the following ODE,

\begin{eqnarray}\label{system 1}
\xi_1^\prime=f_1(\xi_1,\xi_2)  \qquad \xi_2^\prime=f_2(x_1,x_2)
\end{eqnarray}

where,

\begin{eqnarray*}
f_1(x_1,x_2)&=&\frac{-w^2 x_2 + 2(1+w^2 x_1^2)^\frac{3}{2}}{1+w^2(x_1^2+x_2^2)} \, x_2+1\cr
& &\cr
f_2(x_1,x_2)&=&\frac{w^2 x_2 - 2(1+w^2 x_1^2)^\frac{3}{2}}{1+w^2(x_1^2+x_2^2)}\, x_1
\end{eqnarray*}

A direct verification shows that if we define

$$h_w:\bfR{2}\to {\bf R} \com{as} h_w(x_1,x_2)=\frac{x_2}{\sqrt{1+w^2 x_1^2}}+x_1^2+x_2^2$$

then, $h_w$ is a first integral of the ODE for $\xi_1$ and $\xi_2$, i.e. for any solution $\xi_1(s)$ and $\xi_2(s)$ of this system, we have that

$$ h_w(\xi_1(s),\xi_2(s))=M \com{where $M$ is a constant}$$

This completes the proof of the Proposition. \eop

{\bf Proof of Theorem \ref{Flat surfaces}:} If we define the functions $\theta,\, \xi_1$ and $\xi_2$ as in the previous theorem, we have that the equation: Gauss curvature equal zero, that is, the equation $eg-f^2=0$ reduces to $\theta^\prime=-\frac{1}{\xi_2}$. Replacing this equation in the equations $\xi_1^\prime=\theta^\prime \xi_2+1,\quad \xi_2^\prime=-\theta^\prime \xi_1$, we obtain that $\xi_1$ and $\xi_2$ satisfy the following ODE,

\begin{eqnarray}
\xi_1^\prime=0  \qquad \xi_2^\prime=\frac{\xi_1}{\xi_2}
\end{eqnarray}

It follows that $\xi_1(s)=\frac{a}{2}$ for some real number $a$. If $a=0$ then $\xi_2$ is also constant different from zero, and the surface $\phi$  is a cylinder. In the case that $a$ is not zero then $\xi_2=\pm\sqrt{a s+b}$ and $\theta(s)=\mp\frac{2\sqrt{as+b}}{a}$. This completes the proof of the theorem. \eop

As mentioned before, {\bf Theorem \ref{Profile Twizzler vs heart-shape curves}} follows directly from Proposition \ref{ODE for the natural coordinates}. {\bf Lemma \ref{Parametrization heart-shaped}} is a direct computation.

{\bf Proof of Lemma \ref{Length fundamental piece} and Theorem \ref{Moduli space Twizzler}:} Let us start by pointing out that for any $M>-\frac{1}{4}$, the minimum and maximum  distance from the origin to the set $h_w^{-1}(M)$ are $r_1(M) =| \frac{\sqrt{1+ 4 M}-1}{2}|$ and $r_2(M) = \frac{\sqrt{1+ 4 M}+ 1}{2}$ respectively. We can see this by using the method of Lagrange multipliers to find the maximum value and minimum value of the function $R(x_1,x_2) = x_1^2+x_2^2$ subject to the restriction $h_w=M$. We can prove that if $(x_1,x_2)$ is one of these critical points, then $x_1=0$ otherwise the equations

$$\frac{\partial R}{\partial x_1}=\lambda\frac{\partial h_w}{\partial x_1}\com{and} \frac{\partial R}{\partial x_2}=\lambda\frac{\partial h_w}{\partial x_2}$$

has not solution. Once we know that $x_1$ must be zero, we obtain from the equation $h_w=M$ that $x_2$ is
either $-\frac{\sqrt{1+4M}+1}{2}$ or $\frac{\sqrt{1+4M}-1}{2}$. Let us continue with the rest of the proof.
Let us assume that $(x(s),y(s))$ are such that the surface (\ref{surfaces}) has constant mean curvature 1.  Since the curve $(\rho_1,\rho_2)$ defined in Lemma \ref{Parametrization heart-shaped} is regular, we have that

 $$\lambda(u)=
 (\frac{d\rho_1}{du})^2+(\frac{d\rho_2}{du})^2$$

 is a periodic positive function. Likewise, since the functions $f_1(x_1,x_2)$ and $f_2(x_1,x_2)$ only vanish simultaneously at $(x_1,x_2)=(0,-\frac{1}{2})$,  we get,

 $$\mu(u)=f_1^2(\rho_1(u),\rho_2(u))+f_2^2(\rho_1(u),\rho_2(u))$$

 is a positive periodic function. Notice that $\xi_1(s)=0$ and $\xi_2(s)=-\frac{1}{2}$ is the only constant solution of the system (\ref{system 1}). For any other solution, since $h_w$ is a first integral of the system, there exists $M>-\frac{1}{4}$ and a function $\sigma(s)$ such that

 $$\xi_1(s)=\rho_1(\sigma(s))\com{and} \xi_2(s)=\rho_2(\sigma(s))$$

 is a solution of the system (\ref{system 1}). From the equations above we have

 \begin{eqnarray}\label{Eq 1 for the derivative of sigma}
  \xi_1^\prime(s)^2+\xi_2^\prime(s)^2 =\lambda(\sigma(s))\, \sigma^\prime(s)^2
 \end{eqnarray}

 On the other hand we have that,

 $$\xi_1^\prime(s)=f_1(\xi_1(s),\xi_2(s))=f_1(\rho_1(\sigma(s)),\rho_2(s))\com{and} \xi_2^\prime(s)=f_2(\xi_1(s),\xi_2(s))=f_2(\rho_1(\sigma(s)),\rho_2(s))$$

 It follows that $\sigma$ is either strictly increasing or strictly decreasing,  WLOG we can assume that $\sigma$ is strictly increasing.  Therefore we  get that,

 $$\sigma^\prime(s)=\sqrt{\frac{\mu(\sigma(s))}{\lambda(\sigma(s))}}$$

If  $\kappa(u)$ is the inverse of the function $\sigma(s)$, we have that

\begin{eqnarray}\label{derivative of kappa}
\kappa^\prime(u)=\frac{1}{\sigma^\prime(\kappa(u))}=\sqrt{\frac{\lambda(u)}{\mu(u)}}
\end{eqnarray}

If we change from the variable $s$ to the variable $u$, that is, if we consider the functions,

$$\tilde{\theta}(u)=\theta(\kappa(u))\quad \tilde{\xi}_1(u)=\xi_1(\kappa(u))
\quad \tilde{\xi}_2(u)=\xi_2(\kappa(u))\quad
\tilde{x}(u)=x(\kappa(u))\quad \tilde{z}(u)=z(\kappa(u)) $$

It follows from Equation (\ref{derivative of kappa})  and Equation (\ref{derivative theta}), that

$$\tilde{\theta}^\prime(u)=\psi(u)\com{where}
\psi(u)=\frac{-w^2\rho_2(u) + 2(1+w^2\rho_1^2(u))^\frac{3}{2}}{1+w^2(\rho_1^2(u)+\rho_2^2(u))}\, \sqrt{\frac{\lambda(u)}{\mu(u)}}$$

Since $\psi(u)$ is a periodic function with period $2\pi$, it follows by the existence and uniqueness theorem of ODEs that if $\tilde{\theta}(2\pi)=\theta_0$, then

\begin{eqnarray}\label{semi periodic property of theta}
\tilde{\theta}(u+2j\pi)=j\theta_0+\tilde{\theta}(u)\com{for any integer $j$}
\end{eqnarray}

Since $|(x(s),z(s))|=|(\xi_1(s),\xi_2(s))|$, we have that the piece of profile curve

$$C_{f\, p}=C_{fundamental\,  piece}=\{(\tilde{x}(u),\tilde{z}(u)): u\in [0,2\pi]\}$$

also satisfies that $r_1(M)=\min \{|\, q|: q\in C_{f\,  p}\}$ and $r_2(M)=\min \{|\, q|: q\in C_{f\, p}\}$.
Using (\ref{x and z in terms of xi1 and xi2}) and (\ref{semi periodic property of theta}) we get that

\begin{eqnarray}\label{semi periodicity of x and y}
\begin{pmatrix}\tilde{x}(u+2j\pi)\cr \tilde{z}(u+2j\pi)\end{pmatrix}\, =\,
R_{\theta_0}^j \begin{pmatrix}\tilde{x}(u)\cr \tilde{z}(u)\end{pmatrix} \com{where} R_{\theta_0}=\begin{pmatrix}\cos(\theta_0) & -\sin(\theta_0)\cr \sin(\theta_0) & \cos(\theta_0)
\end{pmatrix}
\end{eqnarray}

The equality above implies that the image of the profile curve can be viewed as the orbit of the group $\{R_{\theta_0}^j\}_{j\in \mathbb{Z}}$ acting on $C_{f\, p}$, that is,

\begin{eqnarray}\label{orbit of the fundamental piece}
C=\{(x(t),z(t)):t\in {\bf R}\}=\{R_{\theta_0}^jp:j\in\mathbb{Z}\com{and} p\in C_{fundamental\, piece}\}
\end{eqnarray}

It follows from the equation above that if $\frac{\theta_0}{2 \pi}$ is a rational number, then $C$ is a properly immersed curve and, if  $\frac{\theta_0}{2 \pi}$ is irrational, then $C$ is dense in the annulus $\{(x_1,x_2):r_1(M) \le\sqrt{x_1^2+x_2^2}\le r_2(M)\}$ when $M\ne 0$, or it is dense in the circle of radius 1 when $M=0$. Therefore, Twizzlers with constant mean curvature 1 have the following property: they are properly immersed, or they are dense in the region contained between two concentric cylinders or they are dense in the interior of a cylinder of radius 1. We can prove that a surface corresponding to an irrational value $\frac{\theta_0}{2 \pi}$ is dense by showing that the profile curve is dense, and we can prove that the profile curve is dense by showing that the intersection of this curve with a circle centered at the origin is either the empty set or dense in the circle. The problem of proving this last statement reduces to that of showing that for any irrational number $\iota$ the set $\{\iota-[n\iota]:n\in \mathbb{Z}\}$ is dense in the interval $[0,1]$ which is a known fact. To finish this proof we notice that since the function $(x(s),z(s))$ is parametrize by arc-length, we have that the length of the fundamental piece is $\kappa(2\pi)=\int_0^{2\pi}\sqrt{\frac{\lambda(u)}{\mu(u)}}du$. Also,
since $\tilde{\theta}^{\prime}(u)=\psi(u)$ we have that $\theta_0=\int_0^{2\pi}\psi(u)du$. \eop

{\bf Proof of Theorem \ref{Delaunay's properties}:} Let us assume that $\mathbb{D}(M)$ is parametrized as

$$\phi(s,t)=(x(s),z(x) \sin(t),z(s) \cos(t))$$

where the profile curve $(x(s),z(s))$ is parametrized by arc-length. A direct verification shows that if $\theta(s)$  is a continuous function such that $x^\prime(s)=\cos(\theta(s))$ and $z^\prime(s)=\sin(\theta(s))$, then the mean curvature of
$\mathbb{D}(M)$ is $\frac{1}{2}(\theta^\prime-\frac{\cos(\theta(s))}{z(s)})$. Since the mean curvature of $\mathbb{D}(M)$ is 1, we have the functions $\theta(s)$ and $z(s)$ satisfy the following ordinary differential equation

$$ \theta^\prime=2+\frac{\cos(\theta)}{z}\com{and} z^\prime=\sin(\theta) $$

This ODE has as a first integral the function $h(z,\theta)=z(\cos(\theta)+z)$. Recall that the function $z(s)$ is always positive. Since the minimum of the function $h$ is $-\frac{1}{4}$, it follows that there exists a non-zero constant $k>-\frac{1}{4}$ such that $h(z(s),\theta(s))=k$.  When $k<0$, the level sets of $h(z,\theta)$ are bounded and therefore
$\mathbb{D}(M)$ will represent an undoloid, when $k>0$ the level sets are not bounded and $\mathbb{D}(M)$  will represent a nodoid. In any case, the $z$-values of the level sets of $h(z,\theta)$ are bounded. A direct computation shows that the maximum and minimum of the $z$-values of the level set $h(z,\theta)=k$ are $\frac{1+\sqrt{1+4k}}{2}$ and $|\frac{1-\sqrt{1+4k}}{2}|$. We can prove that $k$ must be equal to $M$ by comparing these critical values of  $z(s)$ with the maximum and the minimum value of the profile curve viewed as the trace of the focus of a conic when it is rolled on a line. A direct computation shows that the Gauss curvature is $-\frac{\theta^\prime\cos(\theta)}{z}$ and since $\theta^\prime=2+\frac{\cos(\theta)}{z}$ then the Gauss curvature reduces to
$-\frac{\cos(\theta)(2z+\cos(\theta))}{z^2}$. Using the Lagrange multiplier method we get that the maximum and the minimun of the Gauss curvature subject to the constrain $h(z,\theta)=k$ are

$$\frac{4\sqrt{1+4 k}}{(1+\sqrt{1+4k})^2} \com{and} -\frac{4\sqrt{1+4 k}}{(-1+\sqrt{1+4k})^2}  $$

respectively. It follows the quotient between the  maximum of the Gauss curvature, and the minimum of the Gauss curvature is $-(\frac{-1+\sqrt{1+4k}}{1+\sqrt{1+4k}})^2$. Since $k=M$, the theorem follows. \eop

{\bf Proof of Theorem \ref{Twizzler's properties}:} Using the same notation as in the proof of Proposition \ref{ODE for the natural coordinates} we get that the Gauss curvature $K$ satisfies

$$K=\frac{eg-f^2}{EG-F^2}=-\frac{w^2(1+\theta^\prime \xi_2)}{(1+w^2\xi_1^2)^2}=-\frac{w^2(1+2\, \xi_2\, \sqrt{1+w^2 \xi_1^2})}{(1+w^2\xi_1^2)(\xi_1^2+\xi_2^2)}$$

Taking  $\rho_1(u)$ and $\rho_2(u)$ as in Lemma \ref{Parametrization heart-shaped} we get the following expression for  the Gauss curvature in terms of the parameter $u$

$$\frac{-4 w^2\sqrt{1+4 M+B \cos^2u}\sin u }{4+w^2+4 A^4 w^4 \cos^4u-2 w^2 \sqrt{1+4 M+B \cos^2u}\sin u+(1+4M)w^2\sin^2u+w^2(8 A^2+B\sin^2u)\cos^2u}$$

A direct computation shows that the derivative of the function $K=K(u)$ is of the form $\cos(u)\, po(u)$ where $po(u)$ is a positive function, therefore the maximum of the Gauss curvature occurs when $u=\frac{3\pi}{2}$ and it is equal to $\frac{2w^2\sqrt{1+4M}}{2+(1+2M+\sqrt{1+4M})w^2}$, and the minimum of the Gauss curvature occurs when $u=\frac{\pi}{2}$ and it is equal to $-\frac{2w^2\sqrt{1+4M}}{2+(1+2M-\sqrt{1+4M})w^2}$. Therefore we get that the quotient of the maximum value of the Gauss curvature and the minimum value of the Gauss curvature is $-\frac{2+(1+2M-\sqrt{1+4M})\, w^2}{2+(1+2M+\sqrt{1+4M})\, w^2}$. The last expression in terms of $v$ transforms into

$$-\frac{1+2M-\sqrt{1+4M}+v-2Mv+ v \, \sqrt{1+4M}}{1+2M+\sqrt{1+4M}+v-2Mv- v \, \sqrt{1+4M}}  $$

A direct verification shows the last expression reduces to $c$ when we replace $v$ by $\frac{\sqrt{1+4M}-1-2M+c(\sqrt{1+4M}+1+2M)}{\sqrt{1+4M}+1-2M+c(\sqrt{1+4M}-1+2M)}$, therefore, all the Twizzlers
$\mathfrak{T}( M,\frac{\sqrt{1+4M}-1-2M+c(\sqrt{1+4M}+1+2M)}{\sqrt{1+4M}+1-2M+c(\sqrt{1+4M}-1+2M)})$ are isometric. This concludes the proof. \eop

{\bf Theorem \ref{special Twizzler's properties}} follows from Theorem \ref{Twizzler's properties} by taking $M=0$.
\vfill
\eject

\section{Some programs using Mathematica}

In this section we will present the file with extension {\it pdf} generated by the software Mathematica for the programs that produce some of the figures in this paper.

\begin{figure}[h]
\centerline{\includegraphics[width=15cm,height=8.84cm]{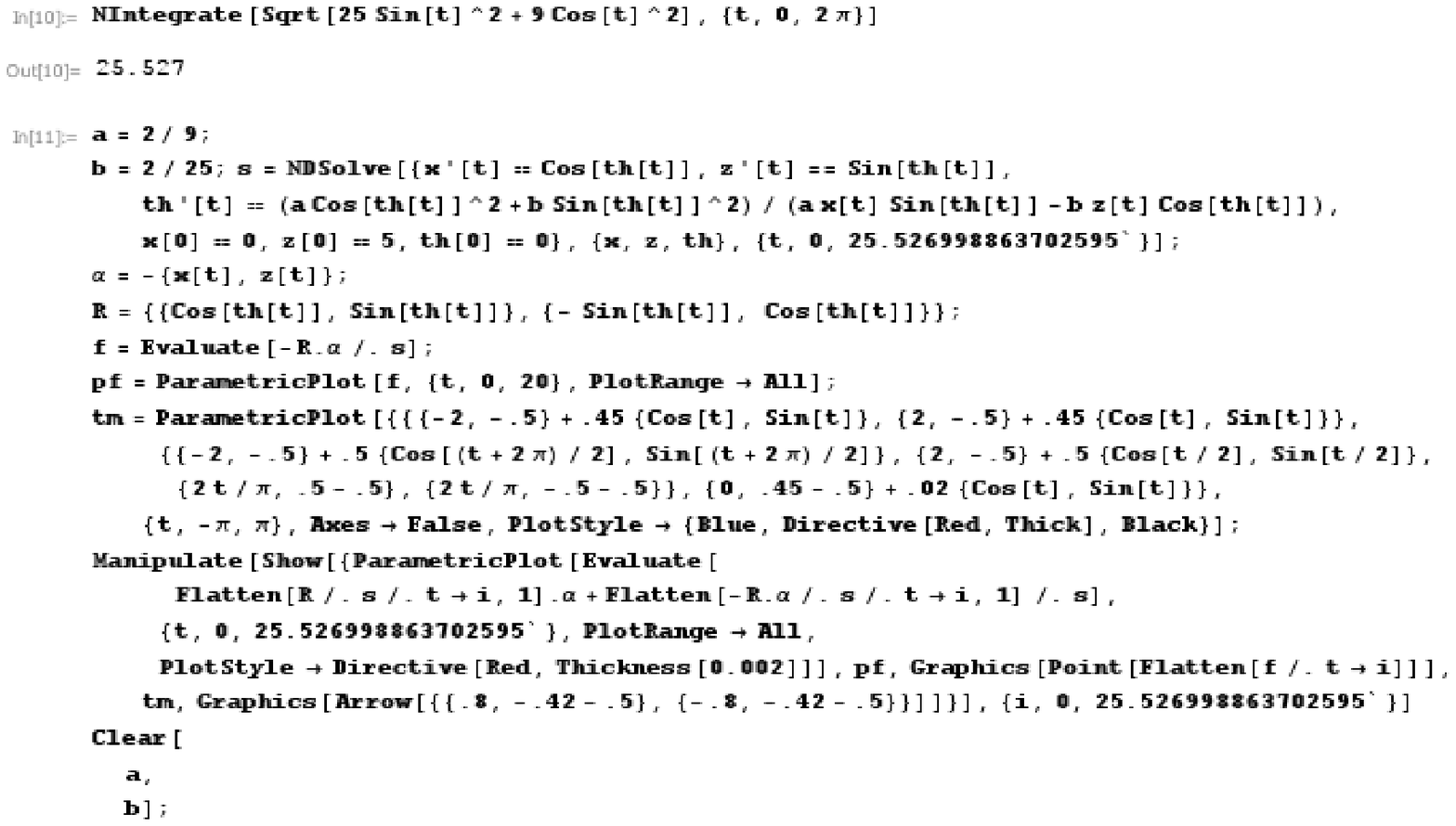}}
\caption{This code produces the video suggested in Figure (1.3)}
\end{figure}

\vfill
\eject

\begin{figure}[h]
\centerline{\includegraphics[width=15cm,height=6.92cm]{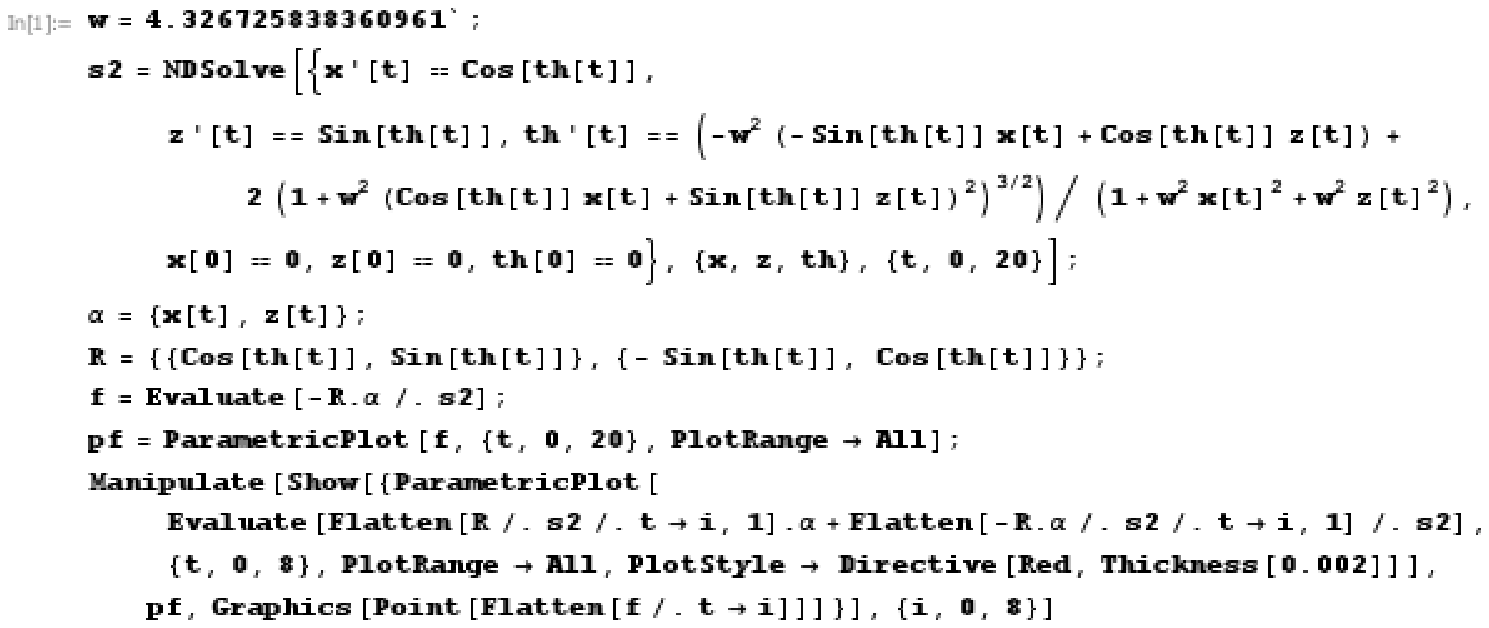}}
\caption{This code produces the video suggested in Figure (1.6)}
\end{figure}

\vskip1cm

\begin{figure}[h]
\centerline{\includegraphics[width=15cm,height=9.95cm]{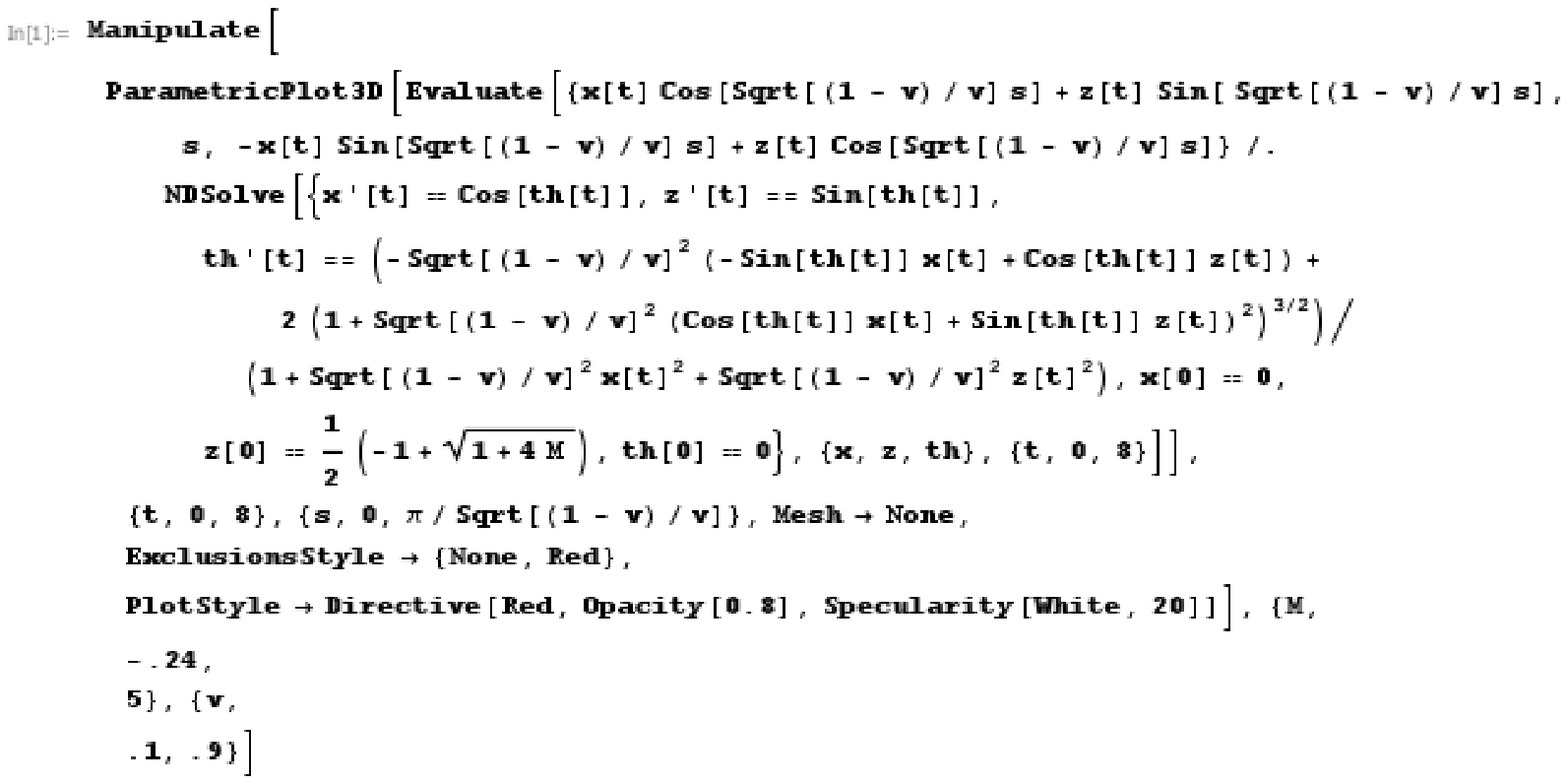}}
\caption{This code produces the video suggested in Figure (1.7)}
\end{figure}


\vfill
\eject

\begin{figure}[h]
\centerline{\includegraphics[width=14cm,height=20cm]{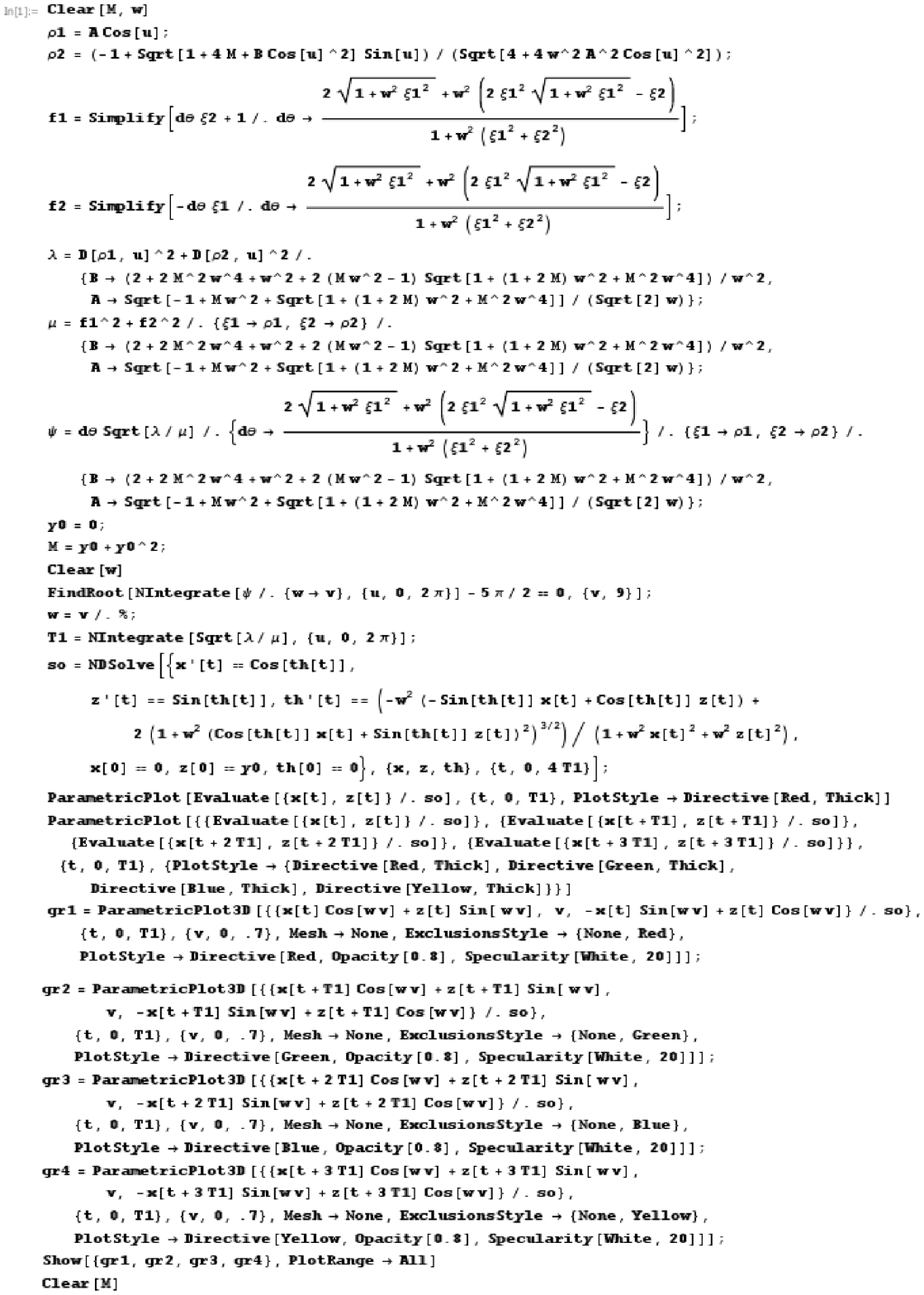}}
\caption{This code produces the picture in Figure (1.8)}
\end{figure}

\vfill
\eject

\begin{figure}[h]
\centerline{\includegraphics[width=15cm,height=9.37cm]{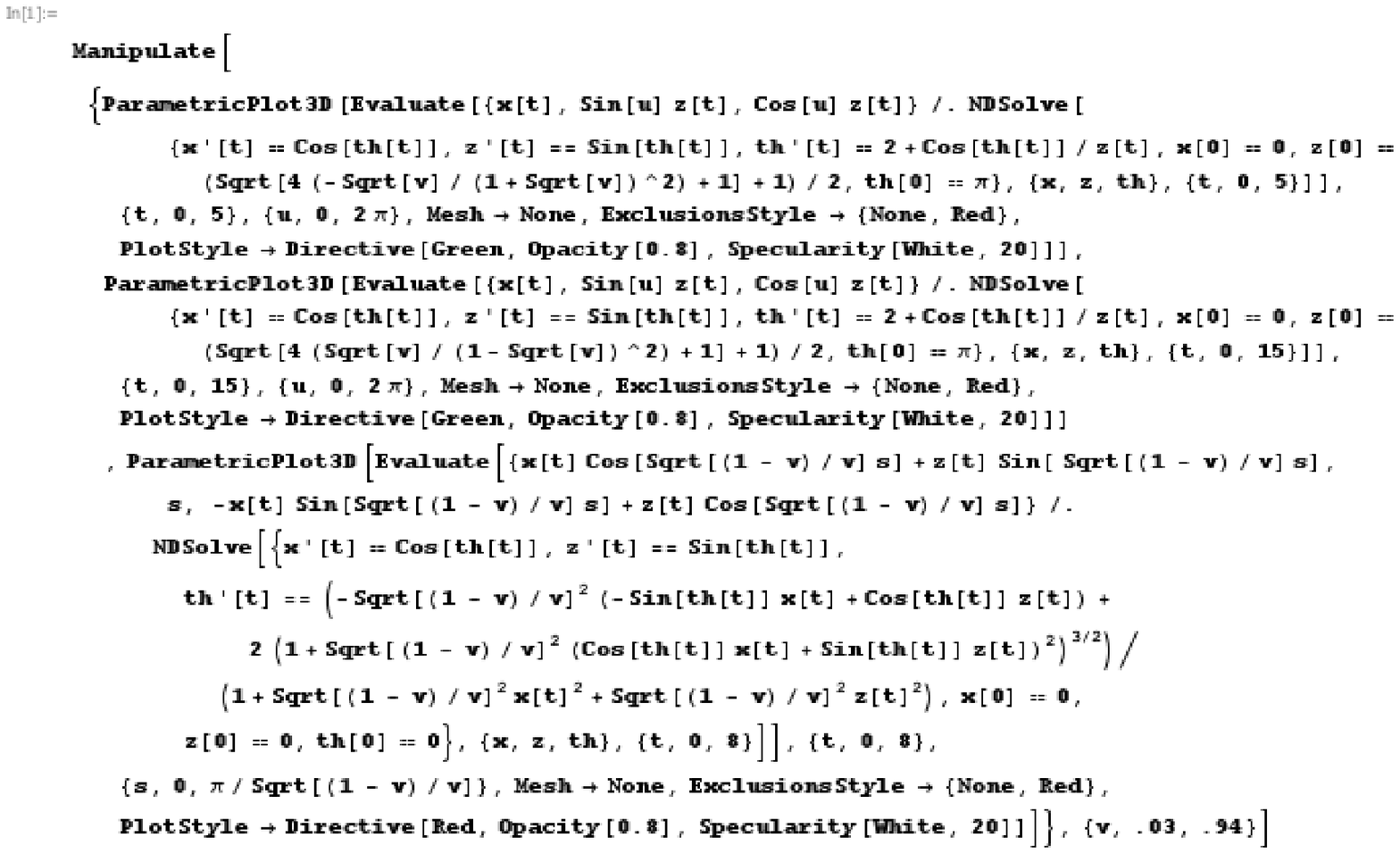}}
\caption{This code produces the video suggested in Figure (1.15)}
\end{figure}

\end{document}